\documentclass[%
 aip,
 amsmath,amssymb,
 reprint,%
]{revtex4-1}

\usepackage{graphicx}
\usepackage{dcolumn}
\usepackage{bm}

\usepackage[utf8]{inputenc}
\usepackage[T1]{fontenc}
\usepackage{mathptmx}
\usepackage{etoolbox}

\newcommand{\RR}{\mathbb{R}}
\newcommand{\XX}{\mathbb{X}}
\newcommand{\YY}{\mathbb{Y}}
\newcommand{\Ddata}{\mathcal{D}}
\newcommand{\bx}{\mathbf{x}}

\usepackage{listings}
\usepackage[capitalize]{cleveref}
\usepackage{amsthm}
\usepackage{float}
\usepackage{graphicx}
\graphicspath{{./Figures/}}

\makeatletter
\def\@email#1#2{%
 \endgroup
 \patchcmd{\titleblock@produce}
  {\frontmatter@RRAPformat}
  {\frontmatter@RRAPformat{\produce@RRAP{*#1\href{mailto:#2}{#2}}}\frontmatter@RRAPformat}
  {}{}
}%
\makeatother
\begin{document}

\preprint{AIP/123-QED}

\title{On the approximation of basins of attraction using deep neural networks}
\author{Joniald Shena*}
\email{jonialdshena@gmail.com}

 \affiliation{National University of Science and Technology ``MISiS'', Moscow, Russia}
\author{Konstantinos Kaloudis}%
\affiliation{ 
Department of Statistics and Actuarial-Financial Mathematics, University of the Aegean, Samos, Greece
}%

\author{Christos Merkatas}
\affiliation{%
Department of Electrical Engineering and Automation, Aalto University, Finland
}%

\author{Miguel A. F. Sanjuán}
\affiliation{Nonlinear Dynamics, Chaos and Complex Systems Group, Departamento de Física, Universidad Rey Juan Carlos, Tulipán s/n, 28933 Móstoles, Madrid, Spain}

\date{\today}

\begin{abstract}
The basin of attraction is the set of initial points that will eventually converge to some attracting set. Its knowledge is important in understanding the dynamical behavior of a given dynamical system of interest. 
In this work, we address the problem of reconstructing the basins of attraction of a multistable system, using only labeled data. To this end, we view this problem as a classification task and use a deep neural network as a classifier for predicting the attractor that corresponds to any given initial condition. Additionally, we provide a method for obtaining an approximation of the basin boundary of the underlying system, using the trained classification model. Finally, we provide evidence relating the complexity of the structure of the basins of attraction with the quality of the obtained reconstructions, via the concept of basin entropy. We demonstrate the application of the proposed method on the Lorenz system in a bistable regime.
\end{abstract}

\maketitle

\begin{quotation}
    In recent years, Machine Learning (ML) algorithms have drawn tremendous attention in various areas of research in nonlinear and complex systems. Interesting applications of ML algorithms include spatiotemporal prediction, prediction of extreme events and chaotic dynamics identification, among others. This research focuses on the problem of reconstructing the basins of attraction of a multistable system, using only a set of initial conditions labeled with their associated attractors. We use these data as a training set for a deep neural network and then exploit the trained network in order to reconstruct the attractors' basins by classifying new points. We illustrate the efficiency of the proposed method in the case of the Lorenz system in a bistable regime and describe the process of additionally obtaining a basin boundary approximation. Finally, we discuss the relation between the drop of the classification accuracy and the increasing complexity of the underlying system, by means of the basin entropy. 
\end{quotation}

\section{\label{sec:introduction}Introduction}
    
    Recent significant advances in the field of Machine Learning (ML), have led to a series of widespread uses of ML based algorithms, ranging from image processing \cite{xie2012image} and speech recognition \cite{hinton2012deep} to mastering chess \cite{campbell2002deep}. Consequently, the general field of nonlinear science has not been unaffected by such advances and the ML approach has been used to address various problems. In particular, after the seminal work of Pathak et al. \cite{pathak2017using, pathak2018model}, where reservoir computing techniques were used in order to reconstruct \cite{pathak2017using} and predict \cite{pathak2018model} complex spatiotemporal dynamics, many important ML-based applications in nonlinear dynamics have been developed. For example, ML has been successfully used for the prediction of extreme events \cite{lellep2020using}, critical transitions \cite{lim2020predicting} and regime changes \cite{brugnago2020predicting}, approximation of the Koopman operator \cite{gulina2021two}, distinguishing regural from chaotic behavior \cite{boulle2020classification}, and identification of chimera states \cite{Ganaie2020}, just to mention a few. One of the main reasons explaining the impressive amount of such contributions, is that the ML framework allows for model-free analysis, without the need of actually knowing the underlying model responsible for the generation of the observed data.
    
    Modeling difficulties might appear in cases where we have at our disposal data generated from an unknown multistable system, that is, a dynamical system having multiple coexisting attractors. Multistable systems arise in many areas of science, as they are connected with interesting phenomena such as pattern formation \cite{mincheva2008multigraph} or tipping points \cite{farazmand2020mitigation}. The ubiquity of multistable systems has attracted the interest of researchers that have proposed a series of ML applications including the optimal design of metamaterials \cite{liu2020machine}, prediction of multistable PDEs with sparse data \cite{chu2021data} and attractor selection \cite{wang2021constrained}. Moreover, recently, Gelbrecht et al. \cite{gelbrecht2021analysis} have proposed a general framework for addressing the predictability of high-dimensional chaotic systems, illustrating their results using a novel Bistable Climate Toy Model. The existence of multiple attractors makes predictability and control of the system more difficult, and it is thus of great importance to devise methods that are able to identify the multiple attractors with high accuracy.
    
    Regarding the model-free analysis of multistable systems, of particular interest is the problem of identifying/reconstructing the structure of the various basins of attraction, using only a set of observed data. In this work, we attempt to solve the problem of the reconstruction of the basins of attraction of high dimensional nonlinear dynamical system by leveraging deep neural networks. Under the proposed framework, reconstructing the basins of attractions can be seen as a classification problem using as training data an appropriately labeled subset of the data, generated directly of the true underlying dynamical system. Using the trained network it is possible to predict the classes of unobserved data, thus obtaining an approximation of the true underlying basins of attraction. In this context, the multiple coexisting attractors constitute the various classes and subsequently, the decision regions (the regions where the data have the same predicted class labels) have the role of the attractor basins. We illustrate the efficiency of the proposed method using the well studied Lorenz system \cite{Lorenz1963} in a bistable regime \cite{cantisan2021transient} and moreover provide evidence that the complexity of the attractor basins has a significant effect on the quality of the obtained reconstruction, using the concept of basin entropy \cite{Daza2016}.   

    Even though the main focus of this work is the reconstruction of the basins of attraction from data, we illustrate that it is also possible to approximate the basin boundary (the manifold that separates the phase space into basins of attraction) of the underlying system, using the trained neural network. To this end, we use proper bounds of the classification probabilities around $0.5$ and keep the associated initial conditions that correspond to these probabilities. Such points are associated with the so-called decision boundaries, which--in our case--are the basin boundaries. We provide evidence that a set consisting of a large enough number of points with classification probability around $0.5$ can be used as a rough approximation of the system basin boundary in the case of the Lorenz flow in a bistable regime. In cases where more accurate approximations of the basin boundary are needed, one could use the obtained set of points lying on the basin boundary and then apply more elaborate techniques \cite{cavoretto2016robust, adamson2003approximating} in order to interpolate a reasonable approximation of the true manifold.

    The article is organized as follows. In Section~\ref{sec:section2}, we analyze the dynamical behavior of the Lorenz system in the regime of interest. In Section~\ref{sec:section3} we briefly review neural networks for classification, along with the details for the generation of the data sets and the numerical simulations. Results for the basins of attraction reconstruction and the basin boundary approximation are presented in Sections~\ref{sec:section3b} and~\ref{sec:section3c}, where we also discuss the relation between the classification accuracy and the basin entropy of the system. In Section~\ref{sec:section4}, we conclude the article with some general remarks and directions for future research.

\section{Model and methods}\label{sec:section2}

    In this section we provide the necessary background for neural network based classification and we review the dynamical behavior of the Lorenz system. 
    
    \subsection{Neural networks} 
        Neural networks are one of the popular approaches to function estimation. It is by now well known, that any neural network with a single hidden layer and an infinite number of hidden neurons can approximate arbitrarily close any function $f:U\subset\RR^m \to \RR^n$ for any compact subset $U\subset \RR^n.$ This property, known as \textit{universal approximation} \cite{cybenko1989approximation} made neural networks attractive to modern machine learning tasks like image classification and estimation of regression functions, to name a few. For our purposes, we aim to use deep neural networks for the reconstruction of the basin of attraction of nonlinear dynamical systems.
        
        For logistic regression, the goal is to learn a mapping $f:\XX\to \YY$ based on available data where $\XX$ is the \textit{feature} space and $\YY$ is the \textit{target} space. Depending on the nature of the space $\YY,$ the problem is called regression or classification. In particular, if $\YY\subseteq \RR^n,\, n\geq 1$ the problem is known as regression while if $\YY$ consists of discrete values, the problem is a classification problem. It is worth noting that in the case where the discrete values imply some kind of \textit{ordering} then the aforementioned problem is known as ordinal regression. However, in this work we are interested in learning the functional relationship between features and values in the set $\{0,1\}.$
        
        In practice, the mapping $f$ is modeled as an artificial neural network (ANN) of $L\geq 1$ layers. In particular, an ANN is a function $f_{\theta}(\cdot)$ parametrized by $\theta=\{W^l,b^l\}_{l=1}^L,$ which represent the collection of weights and biases between the neurons in different layers. Then, based on available data $\Ddata=\{(\bx_i,y_i)\mid 1\leq i \leq N\},$ the aim is to find optimal values for the network parameters in such a way that a cost function is minimized. In particular, for binary classification problems, the objective function is given by \begin{equation}
            L(\theta) = -\frac{1}{N}\sum_{i=1}^N [y_i \log f_\theta(\bx_i) + (1 - y_i)\log(1-f_\theta(\bx_i))] \label{eq:cross-entropy-loss} ,
        \end{equation}
        where $y_i$ is the true value and $f_\theta(\bx_i)$ is the output of the network for the input $\bx_i,$ computed as follows. For each layer $1\leq l \leq L,$ an affine transformation of the input $\bx_i^{l-1}$ is computed recursively as $a^{l} = W^{l}\bx_i^{l-1} + b^l.$ The network applies a nonlinear function $\phi^l,$ known as activation function, to each element of $a^l$ such that $x^l := \phi^l(a^l).$ For a given input $x^0:=\bx$ the output of the network is $f_\theta(\bx):=\phi^l(\bx).$ Typical choices for the activation function $\phi(\cdot)$ are sigmoid functions like the $\tanh(\cdot)$ or the logistic function and, more recently the rectified linear unit (ReLU) \cite{calin2020activation}.
        More details for the derivation of the particular loss function can be found in machine learning references, see Ref.~\onlinecite{bishop2006pattern} for example. 
        
        Optimization of \cref{eq:cross-entropy-loss} is achieved using some gradient based optimization technique \cite{goodfellow2016deep}. The gradients of the above loss function with respect to the parameters of the neural network $\theta$ can be computed via back-propagating the errors through the layers of the network at each step of the optimization procedure, a procedure also known as back-prop. However, most modern deep learning libraries take advantage of \textit{automatic differentiation} (AD) techniques for efficient computation of the gradients of interest. 
        
        After the minimization of the cost function, the optimal set of network parameters $\hat{\theta}:=\{\hat{W}_l,\hat{b}_l\}_{l=1}^{L}$ can be used for prediction of unobserved data. In particular, for a new observation $\mathbf{x}^*$ the predicted class is given from the output of the network at $\mathbf{x}^*$ using the estimated weights. More formally, the output is computed by
        \begin{align}
            \hat{y} = f_{\hat{\theta}}(\mathbf{x}^*).
        \end{align}
        
        For the evaluation of the performance of the classification, several metrics have been proposed in the ML literature. In the following, we will use the accuracy metric defined by
        \begin{align}
            \mathrm{acc} = \frac{\mathrm{\# \text{ of examples correctly classified}}}{\text{\# examples}}.
        \end{align}
        The architecture of the network we are going to use will be specified in the experimental section.
        
    \subsection{Bifurcation diagram and multistability of the Lorenz model}
    
        The Lorenz system \cite{Lorenz1963, layek2015introduction} is a paradigmatic nonlinear dynamical system resulting from the two-dimensional convection occurring in a layer of a fluid of uniform depth, when the temperature difference between the upper and the lower surfaces is kept constant. The mathematical model for this physical system is described from the following system of ordinary differential equations (ODEs):

        \begin{equation}\label{eq:lorenz_odes}
            \begin{aligned}
                \dot{x} &= \sigma \left( y - x \right),  \\ 
                \dot{y} &= r x - y - x z \\ 
                \dot{z} &= x y - \beta z,
            \end{aligned}
        \end{equation}
        with initial condition given by $x(0)=x_0, y(0)=y_0, z(0)=z_0.$
        The constants $\sigma, r$ and $\beta$ are parameters depending on the physical characteristics of the system \cite{Lorenz1963}. 
   
\begin{figure}[htb!]
            \centering
           \includegraphics[width=0.5\textwidth]{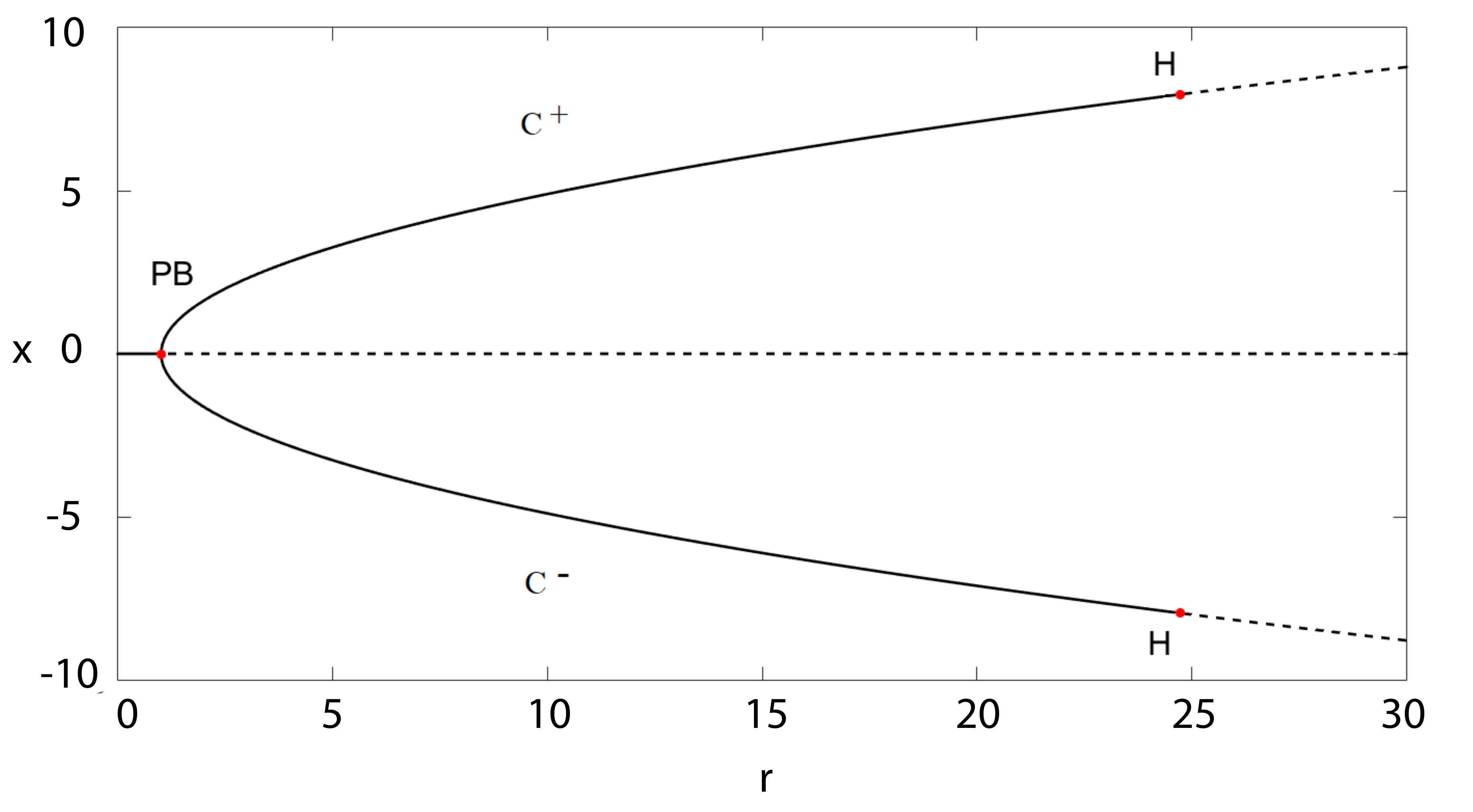}
           \caption{\label{fig:bifurcation} Bifurcation diagram of the Lorenz system. X-coordinate versus the parameter $r$, with $\sigma = 10$ and $\beta = 8/3$. The continuous and dotted lines denote stable and unstable equilibrium respectively. For $r>1$ and after a supercritical pitchfork bifurcation (PB), the system has two stable states, the fixed points attractors $C^{\pm}$. For $r>24.74$ it loses its stability through a subcritical Hopf bifurcation.}
        \end{figure}

        Setting the parameters to the classical values $\sigma = 10$ and $\beta = 8/3$, the dynamics of the system can be more typically explored via a bifurcation diagram for different values of the free parameter $r.$
        Stability analysis of the system is illustrated in Fig.~\ref{fig:bifurcation}. 
        For $0< r < 1,$ the only globally stable attractor is the origin fixed point $(x_0,y_0,z_0)= (0,0,0).$ When $1 \leq r < 24.74$ the system results in one of the two stable fixed points ($C^{-}, C^{+}$) depending on the initial conditions $(x_0,y_0,z_0)$, through a supercritical pitchfork bifurcation (PB). Then, for $r>24.74$ the system passes to chaos through a subcritical Hopf bifurcation.
   
        The Lorenz system is a typical example of a multistable system. Although in this paper our focus is on the region of bistability ($1 \leq r < 24.74$) it should be mentioned that a more detailed stability analysis reveals additional interesting dynamical phenomena. In particular, for $r=24.06$ there is the birth of a chaotic attractor. For $24.06 < r < 24.74$ there is a multistable regime where the chaotic attractors $C^{+}$ and $C^{-}$ coexist \cite{cantisan2021transient} and finally, a pair of unstable limit cycles--called homoclinic orbits--for $13.926<r<24.06$. \cite{Doedel2006}.
   
        In this particular region, except bistability another important phenomenon is transient chaos \cite{lai2011transient}. Transient chaos occurs at $13.926<r<24.06$, where a homoclinic bifurcation takes place and a chaotic saddle is born \cite{Tel2006}. It is related to the stable and unstable manifolds of the chaotic saddle where points exactly on the stable manifold necessarily reach the chaotic saddle and never leave it \cite{cantisan2021transient}. This phenomenon gives rise to fractal regions at the boundaries of the basin of attraction between $C^{+}$ and $C^{-}$.  As we will see below, this behavior is tied with the accuracy of the reconstruction of the basin of attraction.

\section{Experimental setup and results}\label{sec:section3}
    Our aim with the proposed method is twofold. On the one hand, we aim to reconstruct the basin of attraction for the Lorenz map using neural networks. In principle, trajectories that have reached one of the two stable attractors $C^+$ and $C^-$ can be used as training data labeled with the corresponding attractor. On the other hand, we couple the predictive capability of the trained ANN with the basin entropy \citep{Daza2016} of the system. In the following sections, we describe the learning process and elaborate on the connection of the network accuracy with the basin entropy. 
    
    \subsection{Reconstruction of the basin of attraction}\label{sec:section3a}
    
        As mentioned above, having a trajectory generated from the system we can associate each point of the trajectory to one of the coexisting attractors.
        
        In order to obtain the training set, it is sufficient to generate trajectories for the dynamical system starting from distinct initial conditions for time length sufficient enough until the trajectory reaches on the two stable fixed points. In particular, we start by computing the trajectories of the dynamical system for randomly chosen initial conditions $(x_0,y_0,z_0) \sim U\left((-50,50)\times(-50,50)\times(-50,50)\right)$, using the Tsitouras method\cite{tsitouras2011runge} with adaptive time step and absolute and relative tolerance equal to $10^{-6}$. Here, $U(A)$ denotes the uniform density over the set $A.$ For each value of the parameter $1<r<24.74$ for which the system has two stable states, we simulate $100000$ points of the trajectory. The final states of the system are in the attracting set so we can  label this solution with values $1$ and $0$ for each one of the two final states. The integration time for trajectories to arrive to the stable fixed points ranges from $500$ to $1200$, depending on $r.$
    
\begin{figure*}[t]
        \centering
            \includegraphics[width=1\textwidth]{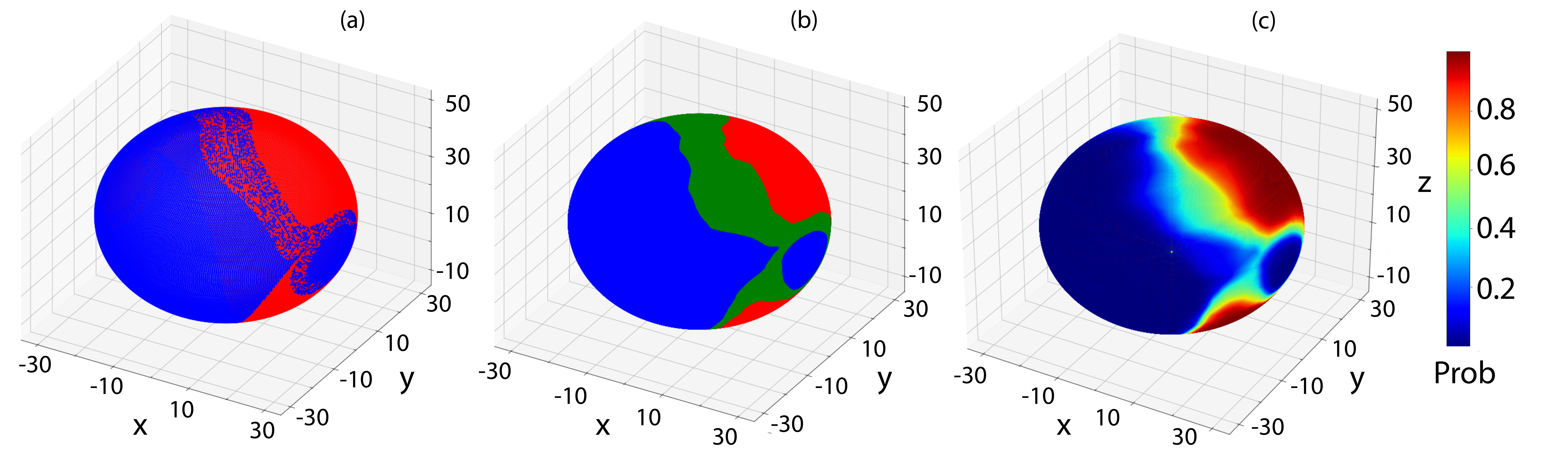}
            \caption{\label{fig:Figure_Final} Basins of attraction for $r = 20$. (a) The sphere of initial conditions where trajectories that end up in the $C^{+}$ and $C^{-}$ have been depicted in red and blue colors, respectively (according to Fig. 2 (a) in Ref.~\onlinecite{cantisan2021transient} ). (b) Classification regions obtained after training the network. Red color refers to a probability greater than 0.75, blue color refers to a probability less than 0.25, and with green color we denote the remaining interval $0.25<Prob<0.75$. (c) More detailed representation of the classification probabilities.}
        \end{figure*}

        Having the trained network, prediction in spherical coordinates can be done with a change of coordinate system. In \cref{fig:Figure_Final}(a), we present the basin of attraction with a sphere with fixed radius $R=30$, centered at $(0, 0, r - 1),$ the midway between $C^+$ and $C^-.$ The connected red colored region corresponds to initial conditions which eventually reach the $C^{+}$ attractor, while in blue color we depict the initial conditions of trajectories that reach the $C^{-}$ attractor. The region where the blue and red dots coexist represents a fractal region. We note here that this is in agreement with the findings in Ref.~\onlinecite{cantisan2021transient}. The reconstruction of the basin of attraction with the proposed method is shown in Fig.~\ref{fig:Figure_Final} (b) where red refers to a probability greater than $0.75,$ blue to a probability less than $0.25,$ and green to the remaining interval $[0.25, 0.75].$ Attractor prediction probabilities are shown in  Fig.~\ref{fig:Figure_Final} (c) in the form of a heat colormap. For the two well defined regions of the basins of attraction \cite{cantisan2021transient}, classification probabilities are close to their extreme values. In contrast, within the mixing area, the neural network is unable to decide the class and results in a probability value close to $0.5.$ We note here, and this will become more clear later on, that there is no substantial benefit in increasing the complexity of the network architecture as the inability to identify the correct region stems from the complexity of the system. 

\begin{figure*}[htb!]
        \centering
            \includegraphics[width=0.8\textwidth]{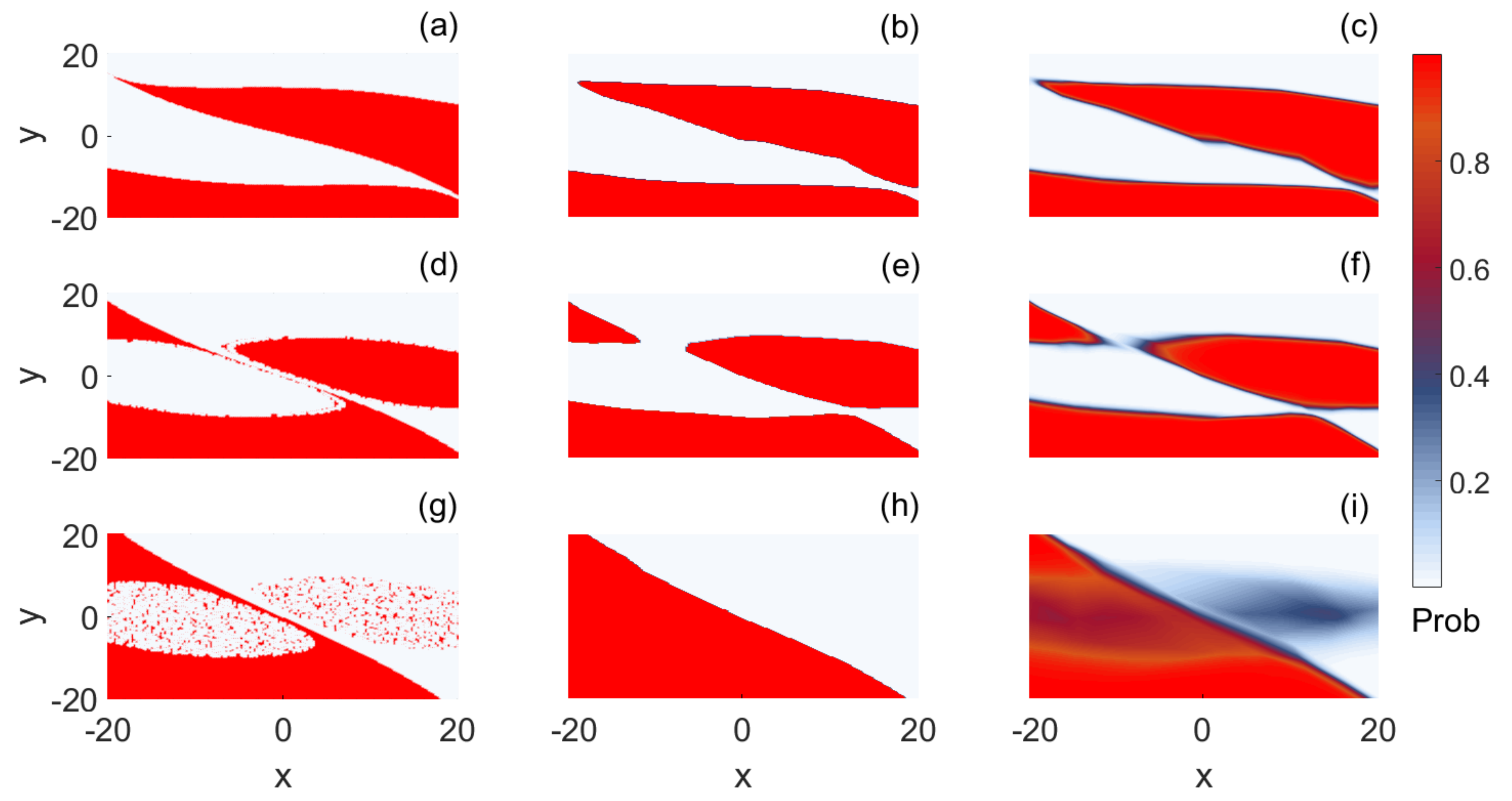}
            \caption{\label{fig:basin_rec_final} A cross section of the basin of attraction for the two equilibria $C^{+}$ and $C^{-}$ (red and gray color, respectively) where $z(t=0) = 5$. The first column represents the actual basin of attraction, the second column refers to the ML estimate, and the third to the probability of that estimate. In the top row ((a), (b) and (c)) $r = 12$, in the middle ((d), (e) and (f)) $r=16$ and in the bottom ((g), (h) and (i)) $r=20$ with ML accuracy achieving $97\%$, $97\%$ and $88\%$, respectively. The color bar illustrates the predicted probability of the ML where equilibria $C^{+}$ and $C^{-}$ corresponds to a probability value equal to $1$ and $0$, respectively.  We have fixed parameters $\sigma = 10$ and $\beta = 8/3$.}
        \end{figure*}

        In Fig.~\ref{fig:basin_rec_final} (a) the basin of attraction is presented for a cross section of the Lorenz system where $z(t=0) = 5$ and $r=12$, for initial conditions in the range $(-20<x<20, -20<y<20)$. The reconstruction of the basin of attraction with the proposed method is shown in Fig.~\ref{fig:basin_rec_final} (b), whilst the probability of that reconstruction is shown in in Fig.~\ref{fig:basin_rec_final} (c). The color bar represents the probability for the particular set of initial conditions to end up in the attracting set $C^{\pm}.$ The ``warmer'' the color the higher the probability. As shown, the trained network can predict very well the basin of attraction ($97\%$ accuracy) with well defined boundaries between the two classes. The probability at the boundaries is close to $0.5$ and this will be crucial later on, for a proper basin boundary approximation.
     
        The same cross section of the basin of attraction for the two equilibria $C^{+}$ and $C^{-}$ is shown in Fig.~\ref{fig:basin_rec_final} (d), (e) and (f), for $r=16$. Now, $r>13.926$ and a homoclinic bifurcation appears leading to the same behavior as in Fig.~\ref{fig:basin_rec_final}. The accuracy of the reconstruction in this case is $97\%$ however, for $r=20$ it decreases significantly. Fractal regions at the boundaries of the basin of attraction between $C^{+}$ and $C^{-}$ take place (Fig.~\ref{fig:basin_rec_final} (g)) leading to low accuracy due to the fact that the two classes cannot be easily separated (Fig.~\ref{fig:basin_rec_final} (h)). Nevertheless, the accuracy of the reconstruction in this case is approximately $88,2\%.$ Experimentally, we have seen that the performance is not so satisfactory in the fractal regions (Fig.~\ref{fig:basin_rec_final} (i)) where the probability of a point to fall in one of the two attractors is close to 0.5.
     
        Next, we test the performance of the proposed method for the whole range of values of the $r$ parameter, where the system displays bistability ($1<r<24.74$). We expect that the accuracy of the reconstruction will decrease as we approach the chaotic region ($r>24.74$). This accuracy can be related to the basin entropy, a new measure for the correct classification of the basin of attraction complexity, for a nonlinear dynamical system \cite{Daza2016}.

    \subsection{Accuracy vs basin entropy}\label{sec:section3b}
    
\begin{figure*}[t]
        \centering
            \includegraphics[width=0.95\textwidth]{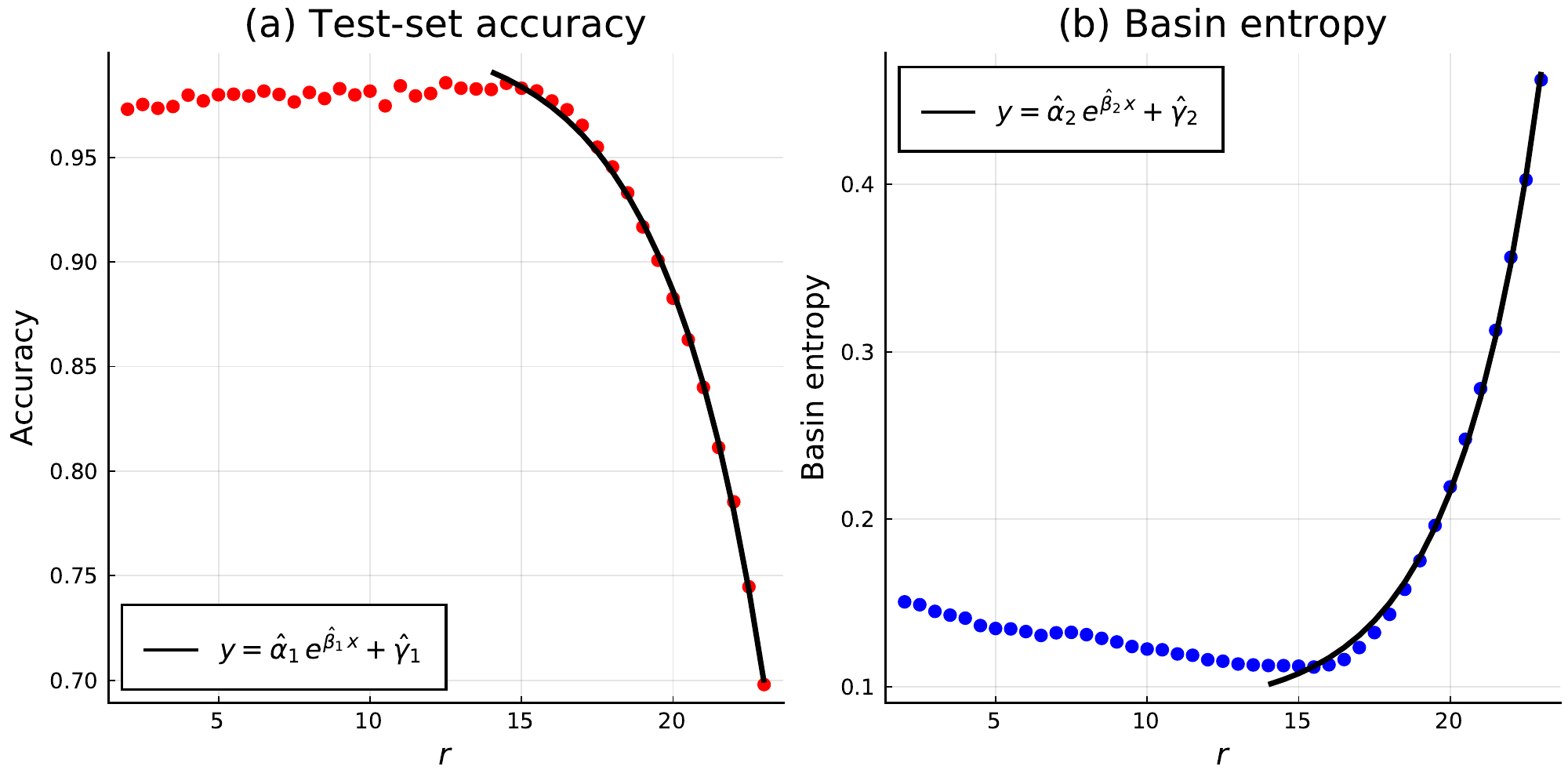}
            \caption{\label{fig:ben_acc} (a) ML classification accuracy measured at the test set and (b) Basin entropy of the Lorenz system, calculated by generating $25$ trajectories per box with a phase space partition at $25$ boxes, superimposed with the fit $f(x)=\alpha\,\exp(\beta\,x)+\gamma$.}
        \end{figure*}

        In order to evaluate the efficiency of the proposed ML classification algorithm with respect to the control parameter $r$ of the system under study, we use the accuracy index. Specifically, for varying values of $r \in \left(2,23\right)$ we train the neural network and then calculate the accuracy index using a data set of $100000$ observations. We split this data set randomly to training and testing sets consisting of $80\%$ and $20\%$ of the observations, respectively. As it has already been mentioned in \cref{sec:section3a}, much care has to be taken regarding the chosen transient times of the generated orbits of the data sets due to the fact that for $r\to 25,$ the transient lifetimes increase substantially. In \cref{fig:ben_acc}(a), we present the obtained accuracy as a function of $r$, showing a clearly decreasing trend as $r$ approaches the chaotic regime. The nonlinear least squares fit of the function $f(x) = \alpha_1\,\exp(\beta_1\,x)+\gamma_1$ indicates an exponential decrease of the accuracy approximately for $r > 14$, with estimated parameters $\left(\hat{\alpha}_1,\hat{\beta}_1,\hat{\gamma}_1\right) = \left( -0.0003, 0.3050, 1.011 \right)$ and corresponding standard errors $\left(3.24\cdot10^{-5}, 0.013,0.004 \right)$. In what follows, we will discuss in more detail the particular emphasis on the values of $r>14$, which are related with the transient chaos onset. 
    
        Interestingly enough, the decreasing trend in the accuracy can be related with the complexity of the associated basins of attraction. In the following, we provide evidence that the structure of the basins of attraction emerges as a natural barrier for the classification accuracy as the dynamical behavior of the system becomes gradually more complex. Similarly as the high Lyapunov exponents are related with the (inevitable) eventual decrease of the prediction accuracy \cite{vallejo2017predictability}.
    
        To this end, we calculate the basin entropy of the system \cite{Daza2016} for $r \in \left(2,23\right)$. The basin entropy among other important features, quantifies the uncertainty of the final state that corresponds to a given value of initial conditions. In order to calculate the basin entropy for each value of $r$, we use a partition of the 3-dimensional phase space in $N=25$ cubes of equal size and generate $25$ orbits for each cube. Then, we compute the basin entropy $S_p = S/N$, where 
        \begin{equation}\label{eq:basin_entropy}
            S = \sum_{i=1}^{N}\sum_{j=1}^{2}p_{i,j}\log\left(1/p_{i,j}\right) ,
        \end{equation}
        where $p_{i,j}$ represents the relative frequency of the $j-$labeled final state occurring from an initial condition of the $i-$box.
    
        The results are presented in \cref{fig:ben_acc}(b), where we observe a clear increasing trend for values of approximately $r>14$. Essentially, this result is a projection of the increasing complexity of the basin structure, thus leading to a relatively higher uncertainty for the system evolution with respect to the choice of the initial condition. Specifically, we perform again a least squares fit of the function $f(x) = \alpha_2\,\exp(\beta_2\,x)+\gamma_2$ that indicates an exponential increase of the $S_b$ for $r > 14$, with estimated parameters $\left(\hat{\alpha}_2,\hat{\beta}_2,\hat{\gamma}_2\right) =  \left(0.0001,0.3570, 0.0859\right)$. The values of the standard errors of the estimated parameters are $\left(6.16\times10^{-5}, 9\times10^{-3}, 3\times10^{-3} \right)$.

\begin{figure*}[t]
        \centering
            \includegraphics[width=0.95\textwidth]{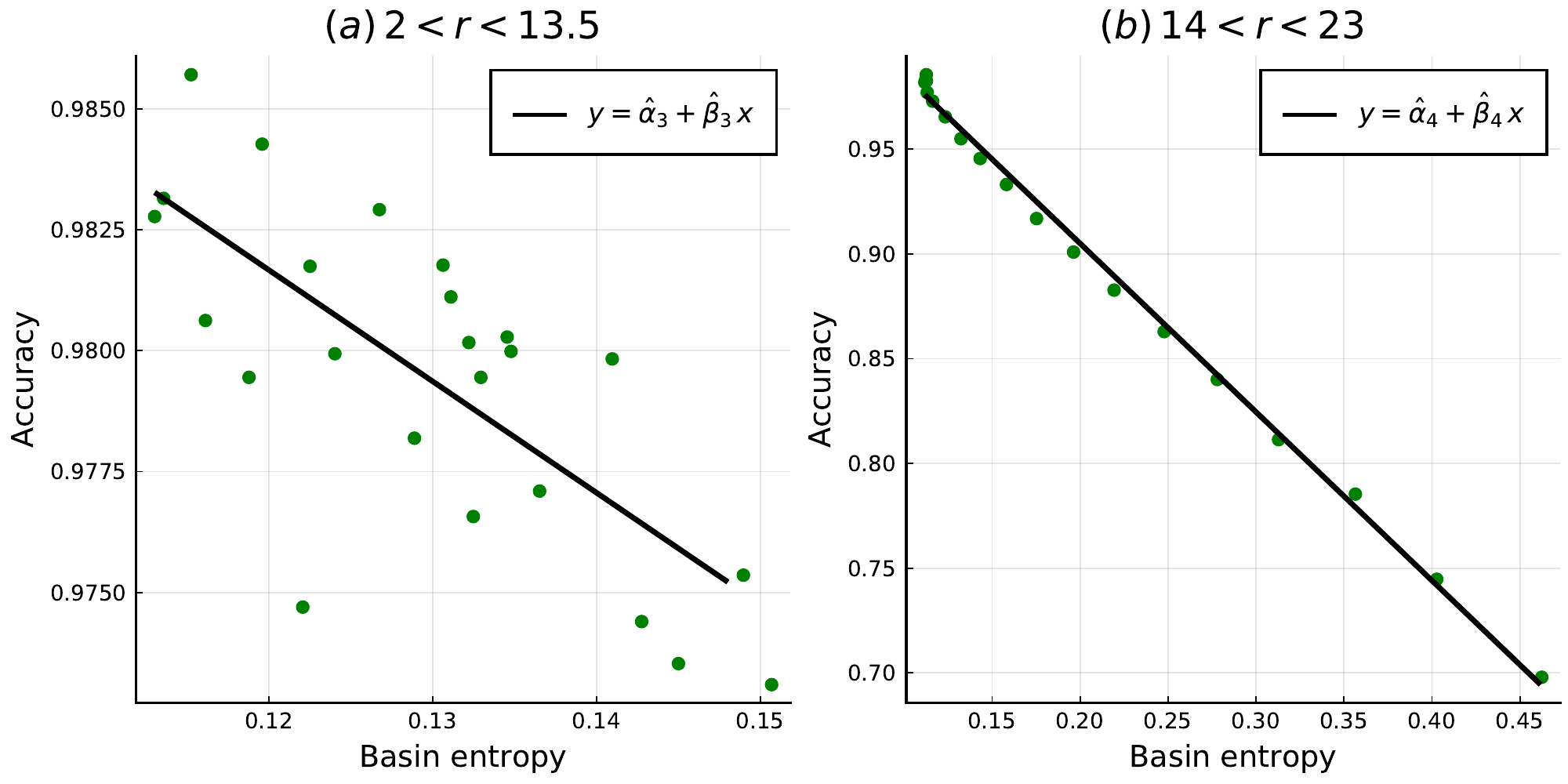}
            \caption{\label{fig:ben_acc2} Scatterplot of accuracy with respect to the basin entropy, for (a) $2<r<13.5$ and (b) $14<r<23$, superimposed with the associated regression lines.}
        \end{figure*}

        The relation between accuracy and basin entropy is further examined in  \cref{fig:ben_acc}(c), with accuracy as a dependent variable. The results of the linear fit $y=\alpha_0+\alpha_1\,x$ indicate a statistically significant linear relation with $\alpha_0=1.085$ and $\alpha_1=-0.85$ (both p-values $< 1\times10^{-10}$). Notably, for low values of the basin entropy (high values of accuracy), there appears to be a more noisy pattern. This fact is partially attributed to the existence of two distinct regions determined by the value of $r$. In particular, in \cref{fig:ben_acc2} we present separately the scatterplots of the accuracy with respect to the basin entropy of \cref{fig:ben_acc} for (a) $2<r<13.5$ and (b) $14<r<23$, superimposed with the associated fitted linear predictors (which are both statistically significant with p-values $1\times10^{-10}$). Although the negative linear correlation is evident in both cases, the first region ($2<r<13.5$) is more noisy and the fitted line is less steep (the values of the estimated slope are $-0.23$ and $-0.8$ for $r<13.5$ and $r>14$, respectively). 
    
         On one hand, the first region is associated with very high values of the classification accuracy, so the pattern of negative correlations is affected by the algorithm's sensitivity to the choice of the specific training and test sets. However, this sensitivity does not justify the high increase of the linear slope. What implies the amplification of the effect that the complexity of the basins has on the classification accuracy. One can deduce that the main reason for this turning point is the occurrence of a homoclinic bifurcation at $r=13.926$ [\onlinecite{cantisan2021transient}]. As described in Ref. [\onlinecite{cantisan2021transient}], the homoclinic bifurcation gives birth to a chaotic saddle, responsible for the transient chaos phenomenon that persists until the emergence of a chaotic attractor at $r=24.06$. Apparently, the birth of the chaotic saddle has an impact on the scaling behavior of the basin entropy, as the final state uncertainty rises quickly after $r\approx14$. It is this uncertainty of the final state of the system that has an amplified effect on the classification accuracy, as compared to the previous regime of $r<13.5$.
    
        Further analyzing this behavior, we can observe that this increase of the complexity of the basins of attraction (equivalently the final state uncertainty) is related to the emergence of fractality at certain regions of the basin boundaries \cite{cantisan2021transient}. In essence, the initial conditions that correspond to long chaotic transients form fractal regions that are reflected at the increase of the basin entropy. This increase implies the drop of the classification accuracy, as it is not possible to classify observations lying in fractal structures due to final state sensitivity \cite{grebogi1983final}. This natural classification accuracy barrier acts similarly to the natural prediction barrier formed by positive Lyapunov exponents \cite{grebogi1983final,vallejo2017predictability}. It is a subject of future research to attempt to construct improved algorithms that will improve (up to a certain level) the classification accuracy, following the attempts to confront the unpredictability issues of large spatiotemporally chaotic systems using echo-state reservoir computing \cite{pathak2018model}.

    \subsection{Basin boundary approximation}\label{sec:section3c}
        As analyzed in the previous section, the complexity of the basins of attractions arises naturally as a barrier to the performance of the classification algorithm, especially at regions where fractality is present. The drop of the classification accuracy however, is related to the associated classification probabilities $p^{\star} \in \left[0,1\right]$.  More specifically, a data point $\mathbf{x}^{\star}$ is classified at a state $y^{\star} \in \left\{0\,(C^{-}),  1\,(C^{+})\right\}$, according to the decision rule
        \begin{align}
            f_{\hat{\theta}}(\mathbf{x^{\star}}) = y^{\star}=\begin{cases} 0,\, p^{\star} < 0.5 \\ 1,\, p^{\star} \geq 0.5 \end{cases}.
        \end{align}

        Thus, the basins of attraction constitute decision regions. Formally, the decision region of the $i$-class label is defined as $DR_i = \left\{\mathbf{x} \in \XX : f(\mathbf{x}) = y_i \right\},\,i = 1,2$ \cite{yan2009studies}. The different decision regions constitute a partition of the state space and are ``separated'' by the decision boundaries. Formally, the intersection of the accumulation points sets of $DR_i$ and $DR_j$ is called decision boundary between the classes $y_i$ and $y_j$ \cite{yan2009studies}. In the context of the present paper, the decision boundaries have an additional meaning, as they represent the basin boundaries which separate the different basins of attraction. The approximation and graphical representation of these manifolds that determine the basins of attraction is a well-known problem that has been addressed in the literature under various approaches, such as meshless approximations \cite{cavoretto2016robust, cavoretto2017graphical}, moving least squares approximants \cite{francomano2018separatrix} and transfer operators \cite{lunsmann2018extended} among others. As a byproduct of the classification algorithm, an approximation of the basin boundary of the underlying unknown system can be obtained through the computation of the decision boundaries.
    
        Regarding the decision boundaries, we note that when the different classes (basins) are separated in a way that has been accurately learned by the classifier, the classification probabilities are close either to $1$ or $0$, except for the points that lie on the boundaries. The results presented in Fig.~\ref{fig:basin_rec_final}, indicate two different cases where the classification probabilities have values around $0.5$, as a sign of high final state (thus, class) uncertainty. The first case consists of basin boundaries that are smooth (see Fig.~\ref{fig:basin_rec_final}(c),(f)), while the second case regards additionally the regions characterized by a fractal structure (see the ``clouds'' appearing in Fig.~\ref{fig:basin_rec_final}(i)). Intuitively, a classification probability around $0.5$ for a point $\mathbf{x}$ that is $\epsilon$-close to a basin boundary, describes an ambiguity of its final state determination. This final state sensitivity \cite{grebogi1983final} means that it may be attracted either to $C^{-}$ or $C^{+}$.
    
        In order to approximate the basin boundary, we generate a dense grid of data and use the trained ANN to estimate the final states. Then, we choose a range around $0.5$ (e.g. $(0.4,0.6)$) and keep all the points that correspond to classification probabilities inside the predefined range. As these points are located on the basin boundaries, the set that contains them can serve as a (rough) approximation of the basin boundary. 
        
\begin{figure}[t]
            \centering
           \includegraphics[width=0.5\textwidth]{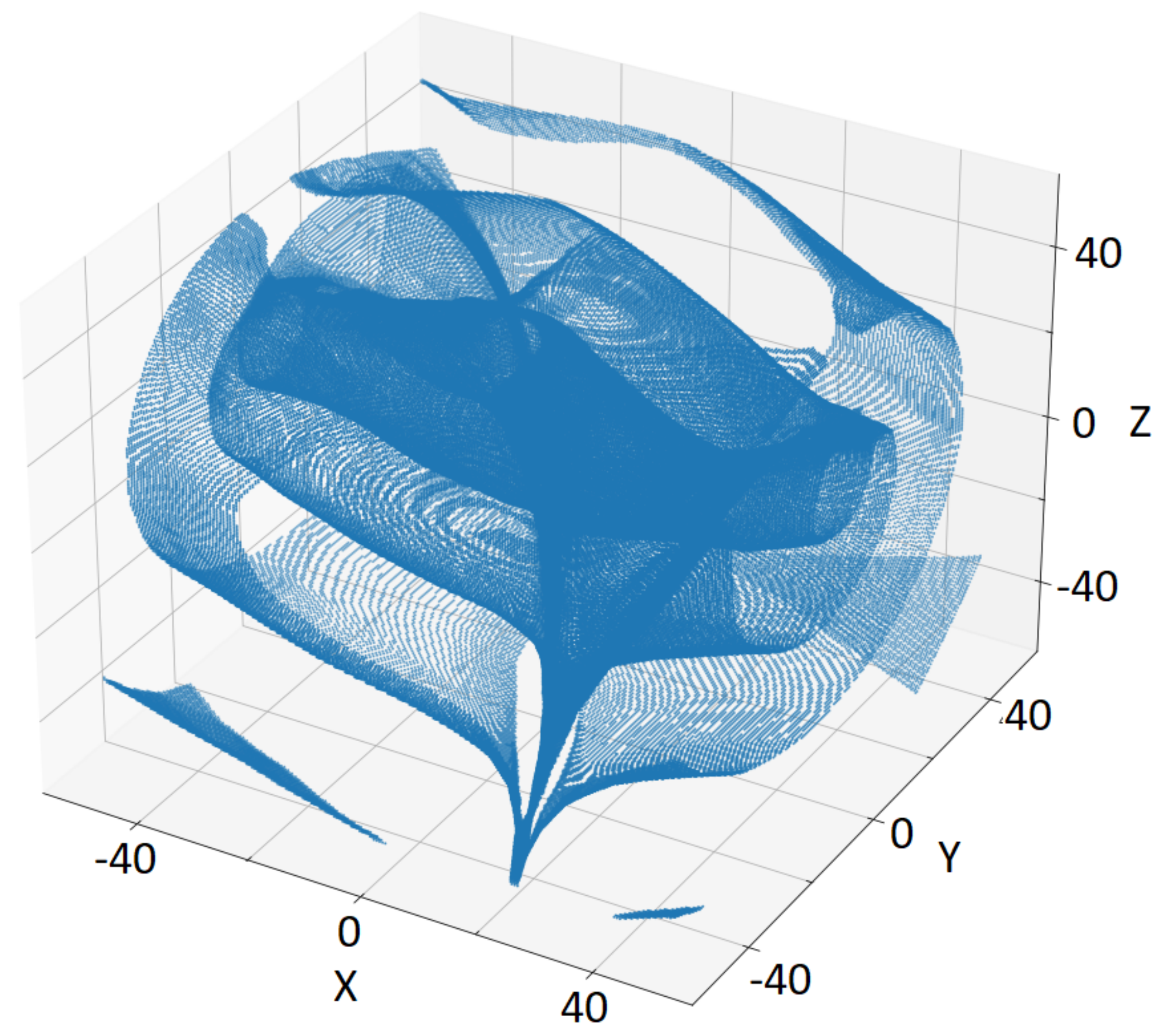}
           \caption{\label{fig:separatrix} The boundary of the basins of coexisting  attractors for the Lorenz system. Parameters have been set to: $r = 12$, $\sigma = 10$ and $\beta = 8/3$.}
        \end{figure}

        In Fig.~\ref{fig:separatrix}, we illustrate the result of the basin boundary approximation in the 3-dimensional state space. The system variables $x, y, z$ range from $-50$ to $50$ and we observe the complicated geometrical characteristics of the basin boundary manifold. The main advantage of the proposed basin boundary approximation -apart from its simplicity- is the numerical efficiency. Having in our disposal the trained classifier, the computational cost for the class prediction of a grid of initial points and the calculation of the classification probability is very low. More specifically, the cost of the assignment of a grid of initial points to their associated class $C^+$ or, $C^-$ reduces to the cost of evaluating the function $f_{\hat{\theta}}(\mathbf{x}^*)$ at the input $\mathbf{x}$.
    
        However, it should be noted that for the basin boundary approximation to be accurate, the range of the control parameter $r$ has to remain between $1$ and $13.92619$, that is,  before the occurrence of the homoclinic bifurcation. More generally, if the system under study exhibits fractal regions at the boundaries of the basins of attraction, it is inevitable that inside these regions the classification probabilities will be around $0.5$ (as a result of the final state sensitivity) and thus, the method will fail to produce a credible approximation. Finally, the proposed method is characterized by wide applicability, as it is a model-free and is able to provide the basin boundary approximation using only a set of labeled data generated by the system of interest.

\section{Conclusions \& future work}\label{sec:section4}
	
	In this work, we have addressed the problem of reconstructing the basins of attraction of a multistable system, using only a set of labeled data, i.e., a set of initial conditions labeled with the associated attractor that will eventually be reached by the orbit starting from the given initial condition. Using this training set, we trained a deep neural network in order to obtain a classification model, which can be used for reconstruction purposes. Specifically, for the reconstruction of the basins of attraction, we have generated a dense grid of initial conditions and obtained the predicted class labels using the classifier (which are the predicted attractors associated with the initial conditions). For the case of the bistable Lorenz system, we have demonstrated the reconstruction of the attractor basins in the 3-dimensional state space with an accuracy close to 97\%. 
	
	Notably, it is straightforward to apply the same classification method in the problem of reconstructing the escape (or exit) basins of a Hamiltonian system such as the H\'{e}non-Heiles \cite{aguirre2001wada}, by labeling each initial condition of the training set according to the associated exit region. Furthermore, one could also apply the same method in cases where the underlying system exhibits more than $2$ stable attractors, resulting in a multilabel classification problem. In such cases, it is possible that more complex network architectures will be needed. If the system under consideration is high dimensional and the number of the coexisting attractors is not known a priori, the additional step of determining them will be necessary, using methods such as Monte Carlo basin bifurcation analysis \cite{gelbrecht2020monte}.
	
	Of course, the quality of the obtained reconstructions depends on the complexity of the true basins of the underlying system. Especially, when we have the emergence of fractal structures in the basin boundaries, the accuracy of the classifier will be lower compared with cases where the boundaries are smooth. In an attempt to investigate this phenomenon in more depth, we have provided evidence that the accuracy of the proposed method is related with the basin entropy, which among other features, quantifies the final state uncertainty of the system. For the system under study, we have found that at the point of the homoclinic bifurcation (thus at the onset of the transient chaos emergence), the increase of the basin entropy is highly correlated with the decrease of the classifier accuracy. Evidently, this relation is further explained by the fractal regions appearing at the basin boundaries, as a result of the emergence of transient chaos (as discussed in detail in Ref. [\onlinecite{cantisan2021transient}]).
	
	Moreover, we highlight the fact that the ``dynamical analogue'' of the decision boundaries in ML terminology are -in this context- the basin boundaries. Utilizing this observation, we illustrated the potential to exploit the trained classification model for the purpose of obtaining an approximation of the basin boundary of the system under study. Although the basin boundary approximation is a difficult task, with the proposed method we showed that a rough approximation of the basin boundary in the 3-dimensional space can easily be obtained by the set of data points associated with classification probabilities around $0.5$. 
	
	Regarding the application of the proposed method for the basin boundary approximation in multilabel classification problems, an additional step has to be taken into consideration. In such cases, instead of a single classification probability we have $k-1$ probabilities for $k$ classes/attractors (as the sum of all $k$ probabilities is equal to $1$). A simple solution to this problem could be, for example, to approximate the basin boundary using the set of points that correspond to a maximum classification probability around $1/k$.
	
	Finally, this work can be seen as a first step in the process of developing more advanced classification oriented ML applications, suitable for more accurate results concerning the model-free analysis of multistable systems, as well as for further studying the effect of the complexity of the basins of attraction on the quality of the obtained approximations.

\begin{acknowledgments}
This work was supported by the Spanish State Research Agency (AEI) and the European Regional Development Fund (FEDER) under Project No. PID2019-105554GB-100 and the Ministry of Education and Science of the Russian Federation in the framework of the Increase Competitiveness Program of NUST “MISiS” (Grant number K4-2018-049).
\end{acknowledgments}

\section*{Data Availability Statement}
The data that support the findings of this study are available from the corresponding author upon reasonable request.

\section*{References}
\bibliography{refs}

\end{document}